\documentclass[10pt,draft]{amsart}
\usepackage{amsmath,amssymb,amsthm}
\begin{document}
\newtheorem{lem}{Lemma}[section]
\newtheorem{prop}{Proposition}[section]
\newtheorem{cor}{Corollary}[section]
\numberwithin{equation}{section}
\newtheorem{thm}{Theorem}[section]
\theoremstyle{remark}
\newtheorem{example}{Example}[section]
\newtheorem*{ack}{Acknowledgment}
\theoremstyle{definition}
\newtheorem{definition}{Definition}[section]
\theoremstyle{remark}
\newtheorem*{notation}{Notation}
\theoremstyle{remark}
\newtheorem{remark}{Remark}[section]
\newenvironment{Abstract}
{\begin{center}\textbf{\footnotesize{Abstract}}%
\end{center} \begin{quote}\begin{footnotesize}}
{\end{footnotesize}\end{quote}\bigskip}
\newenvironment{nome}
{\begin{center}\textbf{{}}%
\end{center} \begin{quote}\end{quote}\bigskip}

\newcommand{\triple}[1]{{|\!|\!|#1|\!|\!|}}
\newcommand{\xx}{\langle x\rangle}
\newcommand{\ep}{\varepsilon}
\newcommand{\al}{\alpha}
\newcommand{\be}{\beta}
\newcommand{\de}{\partial}
\newcommand{\la}{\lambda}
\newcommand{\La}{\Lambda}
\newcommand{\ga}{\gamma}
\newcommand{\del}{\delta}
\newcommand{\Del}{\Delta}
\newcommand{\sig}{\sigma}
\newcommand{\ome}{\omega}
\newcommand{\Ome}{\Omega}
\newcommand{\C}{{\mathbb C}}
\newcommand{\N}{{\mathbb N}}
\newcommand{\Z}{{\mathbb Z}}
\newcommand{\R}{{\mathbb R}}
\newcommand{\Rn}{{\mathbb R}^{n}}
\newcommand{\Rnu}{{\mathbb R}^{n+1}_{+}}
\newcommand{\Cn}{{\mathbb C}^{n}}
\newcommand{\spt}{\,\mathrm{supp}\,}
\newcommand{\Lin}{\mathcal{L}}
\newcommand{\SSS}{\mathcal{S}}
\newcommand{\F}{\mathcal{F}}
\newcommand{\xxi}{\langle\xi\rangle}
\newcommand{\eei}{\langle\eta\rangle}
\newcommand{\xei}{\langle\xi-\eta\rangle}
\newcommand{\yy}{\langle y\rangle}
\newcommand{\dint}{\int\!\!\int}
\newcommand{\hatp}{\widehat\psi}
\renewcommand{\Re}{\;\mathrm{Re}\;}
\renewcommand{\Im}{\;\mathrm{Im}\;}

\title[Scattering for NLS with periodic potential in 1D]
{{ Scattering for small energy solutions of NLS  with periodic
potential  in 1D }}
\author{}

\author{}

\author[Scipio Cuccagna and Nicola Visciglia]{{Scipio Cuccagna\\
DISMI Universit\`a di Modena e Reggio Emilia\\
via Amendola 2, Padiglione Morselli
42100 Reggio Emilia, Italy\\
email: cuccagna.scipio@unimore.it\\
% \vspace{0,2cm}\\
\vspace{0.2cm}
%\author{Nicola Visciglia}
and\\
Nicola Visciglia\\
Dipartimento di Matematica Universit\`a di Pisa\\
Largo B. Pontecorvo 5, 56100 Pisa, Italy\\
email: viscigli@dm.unipi.it\\
tel.: ++39-0502212294, fax: ++39-0502213224}}

%\author[Scipio Cuccagna and Nicola Visciglia]{{Scipio Cuccagna
%and Nicola Visciglia}}

\maketitle

\date{}

\begin{Abstract}
Given $H\equiv -\partial_x^2+V(x)$ with $V:{  \mathbb{R}}\rightarrow
{ \mathbb{R}} $   a smooth    periodic potential, for $\mu \in
\mathbb{R}\backslash \{ 0 \}$ and $p\ge 7$, we prove scattering for
small   solutions to
\begin{equation*}
i \partial_t u +H u=   \mu |u| ^{p-1} u,  \hbox{ } (t, x)\in {
   \mathbb{ {R}} }\times \R, \quad u(0)=u_0\in
H^1({\mathbb{R}}).
\end{equation*}
\end{Abstract}

\section{Introduction}

In this paper, for $\beta:{ \mathbb{R}}^+\rightarrow \mathbb{R}$  a
suitable nonlinearity, we prove scattering of small solutions of
\begin{equation}\label{nls}
i \partial_t u +H u=   \beta(|u|^2)u,  \hbox{ } (t, x)\in \mathbb{
{R}} \times \R, \quad u(0)=u_0\in H^1({\Bbb R})
\end{equation} where $H\equiv -\partial_x^2+V(x)$ with $V(x) $   a smooth real valued
periodic potential. To do this we need to write appropriate
Strichartz  estimates for $H$.
 For every $1\leq p, q\leq \infty$ we consider the
  Birman-Solomjak  spaces
\begin{equation}\label{mixednorm}
l^p({  \mathbb{Z}}, L^q_t[n,n+1])\equiv \left \{f\in L^q_{loc}({\Bbb
R}) \hbox{ s.t. } \{\|f\|_{L^q[n, n+1]} \}_{n\in {\Bbb Z}} \in
l^p({\Bbb Z})\right \},
\end{equation} endowed with the natural norms
\begin{equation*}\begin{aligned} &  \| f \|_{l^p({  \mathbb{Z}}, L^q_t[n,n+1])} ^p\equiv
\sum_{n\in   \mathbb{Z}}   \| f \| _{L^q_t[n,n+1] }^{p}   \hbox{ }
\forall \hbox{ } 1\le p<\infty\hbox{ and } 1\leq q \leq \infty  \\ &
\| f \| _{l^\infty ({ \mathbb{Z}}, L^q_t[n,n+1])} \equiv \sup _{n\in
\mathbb{Z}} \| f \| _{L^q[n,n+1] }.
\end{aligned}
\end{equation*}
 We consider   the Sobolev spaces
\begin{equation}\label{Sobolev}
W^{k,p}({\Bbb R})\equiv \{f\in {\mathcal S}'({\Bbb R})| (1-\partial
_x^2) ^{k/2}f\in L^p(\mathbb{R})\}.
\end{equation}
  For $p=2$ we   set
$H^k(\mathbb{R})\equiv W^{k,2}(\mathbb{R})$. Then we prove:

\begin{thm}\label{thm:nonlin}  Assume $\beta (t)\in C^3(  \mathbb{R},   \mathbb{R}^3)$
with $\beta (0)=\beta '(0)=\beta ''(0)=0 $ and that $V(x)  $ is a
smooth periodic and nonconstant real valued potential. Then there
exists $\epsilon _0>0$
 such that for any initial data $u_0\in H^1({\Bbb
R})$ with $ \| u_0\| _{H^1({\Bbb R})}<\epsilon _0$ problem (\ref{nls})
is globally well--posed. Moreover there exists $C=C(\epsilon _0)>0$
such that
 it is possible to split
$u(t,x)=u_1(t,x)+u_2(t,x)$ so that for any couple $(r,p) $ that satisfies
\begin{equation}\label{condition} 2/{r }+ 1/{p }=1/2 \text{ and }
 (r,p) \in [4, \infty]\times [2, \infty]  ,\end{equation}
 we have
\begin{equation}
\label{ineq:nonlinStricharz}
  \|u_1(t,x) \|    _{\ell ^{\frac{3}{2}r }(
\mathbb{Z}  ,L^{\infty }_{t}([n,n+1],W ^{1,p}({\Bbb R})))} +
\|u_2(t,x) \|    _{ L^{r }_{t}( \mathbb{R},W ^{1,p}({\Bbb
R}))}\leq C \|u_0\|_{H^1({\Bbb R})} .\end{equation}   Furthermore,
there exist  $u_\pm \in H^1({\Bbb R})$ with $ \| u_\pm \|
_{H^1({\Bbb R})}<C\| u_0\| _{H^1({\Bbb R})}$ such that
\begin{equation}
\label{AsFree}\lim _{t\to \pm \infty} \| u(t,x)-e^{-itH} u_\pm \|
_{H^1({\Bbb R})}=0.\end{equation}
\end{thm}

If $V(x)  $   is constant there is a considerable literature on
\eqref{nls}.  A basic tool   are the   Strichartz estimates, see
\cite{cazenave,kt}, which follow, for $\mathcal {V}(t )\equiv
e^{it\partial_x^2}$, from

\begin{equation}\label{dispersivefree}
\| \mathcal {V}(t ) f\|_{L^\infty({\Bbb R})} \leq   C |t|^{-\frac
12} \|f\|_{L^1({\Bbb R})}.
\end{equation}
For any $V(x)$ not constant \eqref{dispersivefree} is not true and
by \cite{cuccagna} we have instead
\begin{equation}\label{dispersive}
\|e^{itH} f\|_{L^\infty({\Bbb R})} \leq C Max\{|t|^{-\frac 12},
\langle t \rangle ^{-\frac 13}\}\|f\|_{L^1({\Bbb R})}.
\end{equation}
\eqref{dispersive} requires a new set of Stricharz estimates for
$e^{itH}$. This is done in the next section. In the subsequent
section we apply the  Stricharz estimates to the nonlinear problem.

In the sequel we shall use the following notations:
$$L^p_x=L^p(\R_x), W^{k,p}_x=W^{k,p}(\R_x), H^s_x=H^s(\R_x).$$

\section{Stricharz estimates}\label{SE}

  For any $r\in [1,
\infty]$ we set $  {r'}= \frac{r}{r-1}.$ By standard arguments it is
possible to prove:
\begin{lem}\label{thmfiniteband}
Let ${\mathcal U}(t):L^2_x\rightarrow L^2_x$ be a uniformly bounded
group in $L^2_x$ such that \ $\|{\mathcal U}(t) f \|_{L^\infty_x}
\leq C_1{\langle t \rangle^{-\frac 13}} \|f\|_{L^1_x}.$ Then there
exists $C>0$ such that for every pair which satisfies
\eqref{condition} we have
\begin{equation}\label{strichartz}
  \|{\mathcal U}(t)f \|_{\ell^{\frac{3}{2}r}( \mathbb{Z},
 L^{\infty}_{t}([n,n+1],L^p_x))} \leq
C \|f\|_{L^2_x}.
\end{equation}
Moreover there is $C>0$ such that for any two pairs $(r_1,p_1)$
  and $(r_2,p_2)$ that satisfy
\eqref{condition}  we have
\begin{equation}\label{duhamel}
\left \| \int _{0}^{t}{\mathcal U}(t-s) F(s) ds\right \| _{\ell
^{\frac 32 r_1}( \mathbb{Z},L^{\infty}_{t}([n,n+1],L ^{p_1}_x))}\end{equation}
$$\leq C \|  F \| _{\ell ^{(\frac 32 r_2)'}(  \mathbb{Z},L^{1}_{t}([n,n+1],L
^{p_2'}_x))}.$$
\end{lem}
Our next step is:
\begin{lem}
 \label{lem:decomposition}
 There exists a projection
 $\pi :L^2_x\to L^2_x$ which commutes with $e^{itH}$
such that the group
  $\mathcal {U}(t )\equiv \pi e^{itH}$ satisfies the hypotheses
  of Lemma \ref{thmfiniteband} and the  group
  $\mathcal {V}(t )\equiv (1-\pi ) e^{itH}$ satisfies the
  estimate \eqref{dispersivefree}.
\end{lem}
\noindent {\bf Proof}. We have $  e^{itH } (x,y)=      K (t,x,y)
 $

 $$ K
(t,x,y)=
  \int   _{ \mathbb{B} }
e^{i(t E  (k) -(x-y)k) }   {m _- ^0(x,k)}{ m  _+ ^0(y,k)} dk
$$
with $e^{\mp i xk}m _\mp ^0(x,k)$ the Bloch functions and $E(k)$ the
band function, see \cite{cuccagna}. By \S 4  \cite{cuccagna}  there
are two characteristic functions $\chi _j(k)$, $j=1,2$ such that
$1=\chi _1(k)+\chi _2(k)$ in $\mathbb{R}$ and such that, if we set

$$ K _j(t,x,y)=
  \int   _{ \mathbb{R} }
e^{i(t E  (k) -(x-y)k) }   {m _- ^0(x,k)}{ m  _+ ^0(y,k)}\chi _j(k)
dk,
$$
then there is a fixed $C>0$ such that $ |K _1(t,x,y)|\le C \langle t
\rangle ^{-\frac{1}{3}}$ and $ |K _2(t,x,y)|\le C   |t|
^{-\frac{1}{2}}$ for all $(t,x,y)\in \mathbb{R}^3$. Notice that
\cite{cuccagna} treats the generic case when all the spectral gaps
of the spectrum  $\sigma (H)$ of  $H$ are nonempty, but the
arguments are the same in the case $\sigma (H)$ has infinitely many
bands with some empty gaps, and much easier if $\sigma (H)$ has
 finitely many bands.

\section{Proof of theorem \ref{thm:nonlin}}\label{nonlinear}

The global well posedness in $H^1_x$ is well know since it
follows from standard theory. Specifically, following a sequence of
arguments in \cite{cazenave} one has:
\begin{itemize}
\item[(1)] if $\| u_0\| _{H^1_x }<\epsilon \le \epsilon _0$ with $\epsilon
_0$ sufficiently small, (\ref{nls}) admits a   solution $u(t)\in
L^\infty _t(\mathbb{R},H^{1}_x) \cap W^{1,\infty }_t (\mathbb{R},
H^{-1}_x)$;

\item[(2)] the above solution is unique;

\item[(3)] the solution $u(t)  $   can be written in the form
$$u(t)=e^{-itH} u_0+ v(t) \text{ with } v(t)=-i\int
_{0}^{t}e^{-i(t-s)H} \beta (|u(s )|^2) u (s )ds.$$

\item[(4)] the above solution is  $u(t)\in
C^0(\mathbb{R}, H^{1}_x) \cap C^1(\mathbb{R}, H^{-1}_x)$ and the
following quantities are conserved:

\begin{align*} &\| u(t)\| _{L^2_x}=\| u_0\| _{L^2_x }, \\&  E(t)
=\int _{\Bbb R} \left ( |\partial _xu(t,x)|^2-V(x) | u(t,x)|^2+2F(|
u(t,x)|^2)\right ) dx=E(0)\end{align*} where $ F(0)=0$ and $\partial
_{\overline{u}}F(| u |^2)=\beta (|u|^2) u;$

\item[(5)] there exists a fixed $C>0$ such that $\| u(t)\|
_{ H^1_x}<C\epsilon $ for all $t\in \mathbb{R}$.
\end{itemize}

Hence we need only to prove the scattering part.  By Lemma
\ref{lem:decomposition} inequality (\ref{ineq:nonlinStricharz}) is
true for some $C=C_0$ for $u$ replaced by $e^{-itH} u_0$. It remains
to show  that (\ref{ineq:nonlinStricharz}) is true with $u$ replaced
in    the left hand side (\ref{ineq:nonlinStricharz})  by the $v$ in
(3). We will show:

\begin{lem}
 \label{lem:continuation}
For $\pi $ the projection in Lemma \ref{lem:decomposition}, let
$v_1(t)=\pi v(t)$ and $v_1(t)=(1-\pi ) v(t)$. Then, for any $D>0$
there are constants $\epsilon _0>0$ and $C(D)$ such that  if $$\|
v_1 (t,x) \|    _{\ell ^{\frac{3}{2}r }( \mathbb{Z}  ,L^{\infty
}_{t}([n,n+1],W ^{1,p}_x))} + \|v_2(t,x) \| _{ L^{r }_{t}(
\mathbb{R},W ^{1,p}_x)}\leq D \| u_0\|_{H^1_x}$$
for all pairs satisfying \eqref{condition}, and if $\|
u_0\|_{H^1_x}<\epsilon <\epsilon _0,$ then
$$\|
v_1 (t,x) \|    _{\ell ^{\frac{3}{2}r }( \mathbb{Z} ,L^{\infty
}_{t}([n,n+1],W ^{1,p}_x))} + \|v_2(t,x) \| _{ L^{r }_{t}(
\mathbb{R},W ^{1,p}_x)}\leq C(D) \epsilon^6 \|u
_0\|_{H^1_x}.$$
\end{lem}

\noindent {\bf Proof.}
 We have
\begin{equation} \label{ineq:NonlinIneq1} \begin{aligned} &  \| v_1  \|    _{\ell ^{\frac{3}{2}r
}( \mathbb{Z},L^{\infty  }_{t}([n,n+1],W ^{1,p }_x))}\lesssim   \|
\beta (|u |^2) u \| _{ L^{ 1}_{t}(\mathbb{R},H ^{1  }_x) } \lesssim
\| |u |^6 u \| _{ L^{ 1}_{t}(\mathbb{R},H ^{1  }_x) } \\ &  \lesssim
\| u \| _{ L^{ \infty }_{t}(\mathbb{R},H ^{1  }_x) } \| u   \| ^6 _{
L^{6}_{t}( \mathbb{R},L ^{\infty }_x) } \le C \|  u_0\| _{H^1_x} \| u
\| ^6 _{ L^{6}_{t}( \mathbb{R},L ^{\infty }_x) }.
\end{aligned}
\end{equation}
Now we split $u=u_1+u_2$ setting $u_1(t)=\pi e^{-itH}u_0+v_1(t)$ and
$u_2(t)=(1-\pi ) e^{-itH}u_0+v_2(t)$. Correspondingly we get by
hypothesis

\begin{equation}\label{ineq:decomposition} \begin{aligned} &
 \| u
\| ^6 _{ L^{6}_{t}( \mathbb{R},L ^{\infty }_x)}\lesssim \| u _1\| ^6
_{ \ell ^6 (\mathbb{Z},L^\infty _t[n, n+1]),L ^{\infty }_x))} +\| u
_2\| ^6 _{ L^{6}_{t}W ^{1,6 }_x }\le C D^6\|u _0\| ^6_{H^1_x}.
\end{aligned}
\end{equation}
By a similar argument
\begin{equation} \label{ineq:NonlinIneq1} \begin{aligned} &
\| v_2  \|    _{  L^{r  }_{t} W ^{1,p }_x }\lesssim   \| \beta (|u
|^2) u \| _{ L^{ 1}_{t}(\mathbb{R},H ^{1  }_x) } \lesssim \| |u |^6
u \| _{ L^{ 1}_{t}(\mathbb{R},H ^{1  }_x) }   \le C D^6\|  u_0\|
_{H^1_x}^7  .
\end{aligned}
\end{equation}
This  yields Lemma \ref{lem:continuation}.

\medskip

The proof of (\ref {AsFree}) is standard and goes as follows.
 $$  e^{itH} u(t)=u_0 -i\int _0^t e^{isH} \beta (|u(s)|^2) u (s)ds $$   and so for
 $t_1<t_2$
 $$  e^{it_2H} u(t_2)-
  e^{it_1H} u(t_1)=-i\int _{t_1}^{t_2}e^{isH}
  \beta (|u(s)|^2) u(s) ds. $$ Then by the proof of Lemma \ref{lem:continuation}

 \begin{equation}  \begin{aligned} & \|e^{it_2H} u(t_2)-  e^{it_1H} u(t_1) \| _{H^1_x}
 \le \| \int _{t_1}^{t_2}e^{isH} \beta (|u(s)|^2) u (s)
 ds \| _{H^1_x}\\& \le \| \beta (|u |^2) u\| _{L^1([t_1,t_2],H^1_x) }\to 0  \text{ for $t_1\to \infty$ and $t_1<t_2$}.\end{aligned}
\end{equation}
Then   $u_+=\lim _{t\to \infty}e^{itH} u(t) $ satisfies the desired
properties. One proves the existence of $u_-=\lim _{t\to
-\infty}e^{itH} u(t) $ similarly.

\end{document}